\documentclass[11pt]{article}
\usepackage{amsfonts,amssymb, amsmath}
\usepackage{graphics}\usepackage{graphicx}
\usepackage{color}
\usepackage{epsfig}
\usepackage{setspace}
\newtheorem{theo}{{\bf{Theorem}}}[section]
\newtheorem{lem}[theo]{{\bf Lemma}}%[section]
\newtheorem{ex}[theo]{{\bf Example}}%[section]
\newtheorem{rem}[theo]{{\bf Remark}}%[section]
\newcommand{\proof}{\noindent {\bf Proof :} }
\newcommand{\tab}{\hspace*{0.3in}}

\newcommand{\qed}{\hfill \rule{1.6mm}{1.6mm}}
\newcommand{\no}{\nonumber}
\newcommand{\noi}{\noindent}
\newcommand{\txt}{\textrm}

\newcommand{\XX}{\mathcal{X}}
\newcommand{\RR}{\mathbb{R}}
\newcommand{\MM}{\mathbb{M}}
\definecolor{co}{rgb}{0.8,0,0.8}
\definecolor{gr}{gray}{0.5}

\newcommand{\va}{\varphi}
\newcommand{\vt}{\vartheta}

\catcode`\@=11

\def\theequation{\@arabic{\c@section}.\@arabic{\c@equation}}
\catcode`\@=12
\title{Marginalization for rare event simulation in switching diffusions \thanks{Support from the European Commission, through the european project iFLY, is gratefully ackowledged}}
\author{Anindya Goswami\thanks{IISER Pune 411008, India, anindya.goswami@gmail.com} \thanks{corresponding author},
 ~~ Fran\c{c}ois Le Gland\thanks{INRIA, Campus de Beaulieu 35042 Rennes C\'{e}dex, France, 	legland@irisa.fr} }
\date{}

\begin{document}
\maketitle
\begin{abstract}
\noi In this paper we use splitting technique to estimate the probability of hitting a rare but critical set by the continuous component of a switching diffusion. Instead of following classical approach we use Wonham filter to achieve multiple goals including reduction of asymptotic variance and exemption from sampling the discrete components.
\end{abstract}

\section{Introduction}
The estimation of rare event probabilities arises in many areas such as insurance, nuclear engineering, air traffic control, communication networks,
etc. Here, we address the problem of estimating the small probability that the continuous component $X(t)$ of a switching diffusion (a special class of a Markov process $(X(t),\theta(t))$ with a hybrid state space) hits a critical region before some final time. We recall that in a switching diffusion, the continuous component $X(t)$ solves a SDE, the coefficients of which may depend on a discrete component $\theta(t)$ that takes a finite number of values (called modes or regimes), while the jump intensity of the discrete component $\theta(t)$ may depend on the continuous component $X(t)$. More precisely $X(t) \in \RR^d$ and $\theta(t) \in \MM= \{1,\ldots, m\}$ where $(X(t), \theta(t))_{t \ge 0}$ is strongly Markov and given by
\begin{eqnarray}\label{9.1}
dX(t)= b(X(t))1_{\theta(t)} dt + dW(t)\\
\label{9.1a} P[\theta(t+h)=j\mid X(t)=x, \theta(t)=i] = \delta_{ij}+ \lambda_{ij}(x)h +o(h)\\
\label{9.1b} (X(0), \theta(0)) \sim \eta
\end{eqnarray}
where $W$ is a Wiener process on $\RR^d$, $b: \RR^d \to \RR^{d\times m}$ a matrix valued measurable function, $\Lambda(x) :=(\lambda_{ij}(x))_{m\times m}$ the rate matrix, $1_i$ the $i$th unit vector in $\RR^m$ and $\eta$ a probability distribution defined on $\RR^d\times \MM$.

The fact that most of the realizations of the underlying process never reach the critical set $B$(say), makes it extremely difficult to estimate the corresponding hitting probability by following classical Monte-Carlo simulation technique. To this end several simulation techniques have been proposed in the literature to estimate the entrance probability into a rare set. We refer \cite{LEC}, \cite{GLA} and references therein for a precise review on these methods.

We adopt the approach that consists in splitting the state space into a sequence of sub-levels, the particle needs to pass before it reaches the rare target. In full generality, splitting methods provide alternative to importance sampling, with the interesting feature that trajectories are simulated under the original model, so that not only the probability of the rare but critical event can be estimated, but also typical critical trajectories
are exhibited. This splitting stage is based on a precise physical description of the evolution of the process between each level leading to the rare set. Following this approach, a Feynman-Kac representation is proposed for the rare event probabilities in \cite{CER} and \cite{CER1}. These papers use and extend the results from \cite{DEL1}, \cite{DEL2}, and \cite{GLA} in order to study the convergence and asymptotic distribution of the related particle system. These works deal with strong Markov processes with general Polish state space. The hybrid process stated in (\ref{9.1})-(\ref{9.1b}) is also a member of this general class of processes. Therefore, automatically \cite{CER} and \cite{CER1} suggests a multi-level splitting approach to estimate the hitting probability of our interest \cite{KRY1},\cite{KRY}. Despite that here we exploit the special structure of hybrid processes to propose a better particle approximation procedure over that naturally follows from existing results.

The approach proposed here is based on marginalization idea, i.e. it relies on using an alternate Markov model $(X(t),\pi(t))$, made of (i) the continuous component $X(t)$ of the switching diffusion and (ii) the conditional probability vector $\pi(t)$ of the discrete component $\theta(t)$ given the past continuous components $X(s)$ for $0\leq s\leq t$. This conditional probability vector is known as the Wonham filter, and jointly with the continuous component, the pair $(X(t),\pi(t))$ is the solution of a SDE. The advantage of this approach is that there is no need to sample the finite set of modes or regimes, which can be a tricky issue, especially if some modes have very small probability. In the context of switching diffusions, there is usually a nominal mode, under which hitting the critical region has a small probability, but there exist also other different ways by which the critical region can be hit. Indeed, there are degraded non-nominal modes under which hitting the critical region is very easy, but switching to these modes has a small probability. The challenge is not only to estimate accurately the overall probability of the rare event, but also to exhibit these different classes of critical trajectories. In particular, we compute the Wonham filters each time by discretizing another set of stochastic differential equations that we derive in the sequel. Some asymptotic results are proved, which show the improvement over a classical splitting method that would sample jointly the continuous and the discrete components of the switching diffusion $(X(t),\theta(t))$.

The rest of this paper is arranged as follows. In Section 2 a system of stochastic differential equations for Wonham filter is derived. Next the multilevel splitting technique is introduced in Section 3. By virtue of splitting method, we construct suitable Markov chains which are later used in Section 4 to obtain Feynman-Kac representations of the hitting probabilities. In Section 5 we select a standard particle approximation scheme and analyze the role of Wonham filter in reduction of asymptotic variance of the particle approximation.
%The variance reduction is illustrated by producing some numerical examples in Section 6.
At the end we add some concluding remarks in Section 6.

\section{Wonham Filter of the Discrete Component}
\begin{lem} \label{9.theo1}Let $\pi = [\pi^1, \ldots, \pi^m]$ belong to $\mathcal{P}(\MM)$, the set of all probabilities on $\MM$. Let $\{(X(t), \theta(t))\}_{t \geq 0}$ be as in (\ref{9.1})- (\ref{9.1a}) and $\{(\tilde{X}(t), \pi(t))\}_{t\geq 0}$ any diffusion process on $\RR^d \times \mathcal{P}(\MM)$, then the following are equivalent.
\begin{itemize}
  \item [(a)] Let ${\cal D}$ denote the set of all smooth functionals on $\RR^d \times \mathcal{P}(\MM)$ which are linear in $\pi$. Let ${\cal L}$ be the differential generator of $\{(\tilde{X}(t), \pi(t))\}_{t\geq 0}$ and ${\cal L}|_{{\cal D}}$ its restriction on ${\cal D}$. Then
\begin{eqnarray} \label{9.L}
	{\cal L}|_{{\cal D}} (\va) = (b(x)\pi)^* \nabla_x \va +(\Lambda(x)^* \pi)^* \nabla _\pi \va + \frac{1}{2} \triangle_x \va +( \sigma(x,\pi){\nabla_x})^*\nabla_\pi \va
\end{eqnarray}
where $\va \in {\cal D}$, $*$ is the transposition, $\nabla_x$ and $\nabla _\pi$ are the column vectors of the first order differential operators w.r.t $x$ and $\pi$ respectively, $\triangle_x$ is the Laplacian w.r.t $x$, $\sigma(x,\pi) = (D(\pi)- \pi\pi^*)b(x)^* $ and $D(\pi)$ is the diagonal matrix with $\pi^i$ as its $i$th diagonal element.
	\item [(b)] For all bounded continuous functional $f$ on $\RR^d\times \mathcal{P}(\MM)$
\begin{eqnarray} \label{9.2}
E[\sum_i \pi^i(t) f(\tilde{X}(t), 1_i) \mid \tilde{X}(s)=x, \pi(s) = \pi]= \sum_i \pi^i E[f(X(t), 1_{\theta(t)}) \mid X(s)=x, \theta(s)=i]
\end{eqnarray}
for every $x \in \RR^d,\pi \in \mathcal{P}(\MM), t \geq s\ge 0$.
\end{itemize}
\end{lem}

\proof Fix $f \in C_b(\RR^d\times \mathcal{P}(\MM))$. Let
\begin{eqnarray}\label{9.8b}
\va(t,x,\pi):= E[\sum_i \pi^i(t) f(\tilde{X}(t), 1_i) \mid \tilde{X}(s)=x, \pi(s) = \pi].
\end{eqnarray}
Using the standard result (see Theorem 3.1.12 \cite{ABH}) of semigroup of operators associated to Markov process, we get that $\va$ solves the following Cauchy problem, provided the solution exists
\begin{eqnarray}
\label{9.8} \frac{\partial}{\partial t} \va (t,x,\pi)&=& 	{\cal L} \va(t,x,\pi)\\
\label{9.8a} \va(s,x,\pi) &=& \sum_i\pi^i f(x, 1_i)
\end{eqnarray}
where $ {\cal L}$ is as in \eqref{9.L}. The fact that the above problem has at most one classical solution follows from the stochastic representation of the solution. Let $\va(t, x, \pi) = \pi^*\psi(t,x)$ be a trial solution, where $\psi: [s,\infty) \times \RR^d \to \RR^m$ is smooth. Then from (\ref{9.8a}), $1_i^*\psi(s,x) = f(x, 1_i)$. Thus from (\ref{9.8}) and (\ref{9.L}) we have
\begin{eqnarray}\label{9.21}
  \pi^* \frac{\partial}{\partial t} \psi - \pi^* b(x)^* \nabla_x (\pi^*\psi) -\pi^* \Lambda(x) \psi -\frac{1}{2} \pi^*\triangle_x \psi
  -(\sigma(x,\pi){\nabla_x})^* \psi% -\frac{1}{2} (\sigma(x,\pi)\sigma(x,\pi)^*{\nabla_\pi})^* \psi
  = 0.
\end{eqnarray}
By substituting the expression of $\sigma(x,\pi)$, we get
\begin{eqnarray*}
 0&=& \pi^*\left( \frac{\partial}{\partial t} \psi - b(x)^* \nabla_x (\pi^*\psi) - \Lambda(x) \psi -\frac{1}{2} \triangle_x \psi \right)
  - \left(D(\pi)b(x)^*\nabla_x\right)^* \psi + \left( \pi\pi^*b(x)^*\nabla_x \right)^* \psi\\
  &=& \pi^*\left( \frac{\partial}{\partial t} \psi - b(x)^* \nabla_x (\pi^*\psi) - \Lambda(x) \psi -\frac{1}{2} \triangle_x \psi\right)
  - \pi^*D\left( b(x)^*\nabla_x\right) \psi + \pi^*b(x)^* \nabla_x (\pi^*\psi)\\
 &=& \pi^*\left( \frac{\partial}{\partial t} \psi - \Lambda(x) \psi -\frac{1}{2} \triangle_x \psi - D\left( b(x)^*\nabla_x\right) \psi\right)
  \end{eqnarray*}
for all $\pi$. Hence we have
\begin{eqnarray}\label{9.6}
\frac{\partial}{\partial t} 1_i^*\psi - 1_i^*b(x)^*\nabla_x (1_i^*\psi) -\frac{1}{2} \triangle_x (1_i^*\psi) - 1_i^*\Lambda(x) \psi &=& 0 ~~ t>s\\
\label{9.6a} 1_i^* \psi(s,x) &=& f(x, 1_i)
\end{eqnarray}
for all $i =1, \ldots, m$. It turns out that the Cauchy problem of the parabolic system (\ref{9.6})-(\ref{9.6a}) has a classical solution. Again, using the differential generator of semigroup associated with the Markov process $(X(t), \theta(t))_{t\ge 0}$ (which is given by (\ref{9.1})-(\ref{9.1a})), we conclude that the solution of (\ref{9.6})-(\ref{9.6a}) has the following stochastic representation
\begin{eqnarray}\label{9.7}
\no 1_i^*\psi(t,x) = E[f(X(t),1_{\theta(t)}) \mid X(s)=x, \theta(s)=i].
\end{eqnarray}
Therefore, (\ref{9.8})-(\ref{9.8a}) has a unique classical solution given by $\va(t,x,\pi) = \sum_i \pi^i E[f(X(t),1_{\theta(t)}) \mid X(s)=x, \theta(s)=i]$. Hence, (a) implies (b).

Let $\{S_t\}_{t\ge s}$ denote the semigroup associated to the Markov process $\{(X(t), \pi(t))\}_{t \geq s}$, i.e.,
\begin{equation*}
S_t f(x,\pi) := E[f(\tilde{X}(t), \pi(t) ) \mid \tilde{X}(s)=x, \pi(s) = \pi].
\end{equation*}
Thus ${\cal L}$ is the differential generator of $\{(X(t), \pi(t))\}_{t \geq s}$ if and only if
\begin{equation*}
\frac{\partial}{\partial t}S_t f = {\cal L} S_t f \tab \forall f \in C_b(\RR^d\times \mathcal{P}(\MM)).
\end{equation*}
Again since $S_t$ is a linear operator for every $t$, we can conclude that ${\cal L}|_{\cal D}$ is the restriction of generator on ${\cal D}$ if and only if
\begin{equation*}
\frac{\partial}{\partial t}S_t f = {\cal L}|_{\cal D} S_t f \tab \forall f \in {\cal D}.
\end{equation*}
From \eqref{9.2} we get for every $f \in {\cal D}$
\begin{eqnarray*}
% \nonumber to remove numbering (before each equation)
  S_t f(x,\pi)&=& E[\sum_i \pi^i(t) f(\tilde{X}(t), 1_i) \mid \tilde{X}(s)=x, \pi(s) = \pi]\\
  &=& \sum_i \pi^i E[f(X(t), 1_{\theta(t)}) \mid X(s)=x, \theta(s)=i]\\
  &=& \pi^* \psi(t,x) \tab \textrm{ (say).}
\end{eqnarray*}
where $\psi$ can clearly be represented as the solution of (\ref{9.6})-(\ref{9.6a}). Therefore, $\psi$ also satisfies \eqref{9.21}. Thus for every $f \in {\cal D}$
\begin{eqnarray*}
% \nonumber to remove numbering (before each equation)
 {\cal L}|_{\cal D} S_t f (x,\pi)&=& \pi^* \frac{\partial}{\partial t} \psi(t,x)\\
  &=&\pi^* b(x)^* \nabla_x (\pi^*\psi(t,x)) +\pi^* \Lambda(x) \psi (t,x) +\frac{1}{2} \pi^*\triangle_x \psi (t,x) +(\sigma(x,\pi){\nabla_x})^* \psi (t,x)\\
  &=& \Big((b(x)\pi)^* \nabla_x +(\Lambda(x)^* \pi)^* \nabla _\pi + \frac{1}{2} \triangle_x +( \sigma(x,\pi){\nabla_x})^*\nabla_\pi \Big) \pi^*\psi(t,x) \\
  &=& \Big((b(x)\pi)^* \nabla_x +(\Lambda(x)^* \pi)^* \nabla _\pi + \frac{1}{2} \triangle_x +( \sigma(x,\pi){\nabla_x})^*\nabla_\pi \Big) S_t f(x,\pi)
\end{eqnarray*}
for every $(x,\pi) \in \RR^d\times \mathcal{P}(\MM)$. Hence (b) implies (a). \qed

\begin{theo} \label{9.rem1} Let
\begin{equation}\label{9.7a}
\pi^i(t) := P(\theta(t) =i\mid {\cal F}^X(t)), ~ i=1, \ldots,m
\end{equation}
where $(X(t), \theta(t))_{t\ge 0}$ is as in \eqref{9.1}-\eqref{9.1b}. Then the diffusion process $(X(t), \pi(t))$ is given by
\begin{eqnarray}\label{9.10}
d X(t) &=& b(X(t)) \pi(t) dt + d \tilde{W}(t)\\
\label{9.10a} d\pi(t) &=& \Lambda(X(t))^* \pi(t) dt + (D(\pi(t))- \pi(t)\pi^*(t))b(X(t))^* d \tilde{W}(t)\\
\label{9.10b}\pi^i(0) &=& \frac{d \eta(\cdot, i)}{d \sum_i \eta(\cdot, i)}(X(0))\\
\label{9.10c}X(0) &\sim& \sum_i \eta(\cdot, i)
\end{eqnarray}
where $\tilde{W}$ is a Weiner process on $\RR^d$.
\end{theo}
\proof
Consider the process $(X(t), \pi(t))$, given by \eqref{9.1}-\eqref{9.1b} and \eqref{9.7a}. Then
\begin{eqnarray*}
% \nonumber to remove numbering (before each equation)
\lefteqn{ E\left[\sum_i \pi^i(t) f(X(t), 1_i) \mid X(s), \pi(s)\right]}\\
  &=& E\left[E\left[\sum_i P(\theta(t) =i\mid {\cal F}^X(t)) f(X(t), 1_i) \mid X(s), \theta(s)\right]\mid {\cal F}^X(s)\right] \\
  &=& \sum_j P(\theta(s) =j\mid {\cal F}^X(s)) E\left[E\left[ f(X(t), 1_{\theta(t)})\mid {\cal F}^X(t)\right] \mid X(s), \theta(s)=j\right] \\
  &=& \sum_j \pi^j(s) E\left[ f(X(t), 1_{\theta(t)})\mid X(s), \theta(s)=j\right].
\end{eqnarray*}
Hence, $\{(X(t), \pi(t))\}_{t\geq 0}$ satisfies \eqref{9.2}. Therefore from Lemma \ref{9.theo1}, its differential generator is given by \eqref{9.L}. Thus the diffusion process $(X(t), \pi(t))$ can be represented by the stochastic differential equations \eqref{9.10}-\eqref{9.10a}. Again from the definition of $X(0)$ and Wonham filter $\pi(0)$, we have \eqref{9.10b}-\eqref{9.10c}.\qed

The process $\{\pi(t)\}_{t\geq 0}$ as in \eqref{9.7a} is called \emph{Wonham filter}. We cite \cite{G} for a detailed study of a special case of \eqref{9.10}-\eqref{9.10a}. Some comments on Wonham filter is given below.

\begin{rem}
Since $\pi(t)$ is a probability vector for all $t\ge 0$, we must have
\begin{eqnarray}\label{9.18}
\hat{1}^* d\pi(t) =0
\end{eqnarray}
for $t \ge 0$ where $\hat{1}:= \sum_{i\in \MM} 1_i$. Therefore it is relevant to check the above necessary condition (\ref{9.18}) from the stochastic differential equation (\ref{9.10a}). Since, $\hat{1}^* D(\pi) = \pi^*$ and $\hat{1}^* \pi= 1$, we have $\hat{1}^*(D(\pi(t))- \pi(t)\pi^*(t))=0$. Again, $\hat{1}^* \Lambda(X(t))^* =0$ follows from the fact that the row sum of the rate matrix $\Lambda(x)$ vanishes for all rows. Thus (\ref{9.18}) is satisfied.
\end{rem}

\begin{rem}\label{9.rem2}
\noi Before moving to the next section, we wish to emphasize on a certain inherent deterministic nature of the Wonham filter in this remark. We define the following two nonlinear smooth functions\\
$F: \RR^d \times \mathcal{P}(\MM) \to \RR^m$ and $ G: \RR^d \times \mathcal{P}(\MM) \to \RR^{m\times d}$ given by $G(x,\pi) :=(D(\pi)- \pi\pi^*)b(x)^*$ and $F(x,\pi) := \Lambda(x)^* \pi -G(x,\pi)b(x)\pi$. Using (\ref{9.10})-(\ref{9.10a}) we obtain
\begin{equation}\label{9.18a}
d\pi(t) = F(X(t), \pi(t)) dt + G(X(t), \pi(t)) dX(t).
\end{equation}
Given $y \in \RR^d, i\in \MM, y^t_{s} \in C([s,t]; \RR^d), y^t_{s}(s)= y$ and $\pi \in \mathcal{P}(\MM)$, define $\pi^s_t(y, i, y^t_{s}) \in \mathcal{P}(\MM)$ by
\begin{eqnarray*}
\pi^s_t(y, i, y^t_{s})(j):= P(\theta(t)=j \mid X(s)=y, \theta(s)=i, X(u)= y^t_{s}(u), u \in [s,t]), ~ j \in \MM
\end{eqnarray*}
and define $(\tilde{\pi}^s_u(y, \pi, y^t_{s}))_{s\le u \le t}$ as the solution of
\begin{eqnarray*}
d \tilde{\pi}(u) &=& F(y^t_{s}(u), \tilde{\pi}(u))du + G(y^t_{s}(u), \tilde{\pi}(u))dy^t_{s}(u)\\
\tilde{\pi}(s) &=& \pi.
\end{eqnarray*}
Now using \eqref{9.18a}, \eqref{9.2} and \eqref{9.7a} we get
\begin{eqnarray*}\label{9.5}
\lefteqn{\tilde{\pi}^s_t(y, \pi, y^t_{s})(j)}\\
&=& E[\pi^j(t) \mid X(s)=y, \pi(s)=\pi, X(u)= y^t_{s}(u), u \in [s,t]]\\
&=& \sum_{i \in \MM}\pi^i E[P[\theta(t)=j \mid \mathcal{F}^X_t]\mid X(s)=y, \theta(s)=i, X(u)= y^t_{s}(u), u \in [s,t]]\\
&=& \sum_{i \in \MM}\pi^i P[\theta(t)=j \mid X(s)=y, \theta(s)=i, X(u)= y^t_{s}(u), u \in [s,t]]\\
&=& \sum_{i \in \MM}\pi^i \pi^s_t(y, i, y^t_{s}).
\end{eqnarray*}
\end{rem}

\section{Multilevel Splitting and Associated Markov Chains}
The objective of this paper is to reduce the variance in estimating the probability that the continuous component $X(t)$ of the hybrid process $\{(X(t), \theta(t))\}_{t \ge 0}$ hits a critical subregion $B\subset \RR^d$. We adopt the multilevel splitting approach, thus the intermediate subregions are introduced in the following way
$$
B=B_n \subset B_{n-1} \subset \cdots \subset B_0 = \RR^d. $$
Let $T_k$, $k=1, \ldots, n$ be the associated hitting times, i.e.;
\begin{eqnarray}\label{9.12}
\no T_k:= \inf\{ t\ge 0 : X(t) \in B_k\}
\end{eqnarray}
and hence $T_k$s are non decreasing $\mathcal{F}^X$-stopping times. Let $S_k:= T_k \wedge T$ be the truncated stopping times with $S_0 \equiv S_{-1} \equiv 0$. We denote $\XX_k: [S_{k-1}, S_k] \to \RR^d$, defined as $\XX_k(t) = X(t)$ for $S_{k-1}\le t \le S_k$ and $k=0,1, \ldots, n$. Thus the random variables $\XX_k$ take values in ${\cal E} := \{y^t_{s} \in C([s,t]; \RR^d) : 0\le s \le t\}$ for $k=0,1, \ldots, n$. The set ${\cal E}$ can be embedded into $\RR^2 \times C([0,T]; \RR^d)$ injectively in the following way
$$ y^t_{s} \mapsto \left(s,t, I_{[0,s)}(\cdot)y^t_{s}(s) + I_{[s,t]}(\cdot)y^t_{s}(\cdot)+ I_{(t,T]}(\cdot)y^t_{s}(t)\right)
$$
and we endow ${\cal E}$ with the $\sigma$-algebra that is inherited from the Borel $\sigma$-algebra of $\RR^2 \times C([0,\infty); \RR^d)$.

We define the sequences $\bar{\theta}_k := \theta(S_k)$ and $\bar{\pi}_k := \pi(S_k)$ where $\theta(t)$ and $\pi(t)$ are given by the equations (\ref{9.1})-(\ref{9.1b}) and (\ref{9.10})-(\ref{9.10c}) respectively.

\begin{theo}\label{9.theo2}
$\xi_k:= (\XX_k,\bar{\pi}_k)$ and $\vt_k:= (\XX_k,\bar{\theta}_k)$ are ${\cal E}\times{\cal P}(\MM)$ and ${\cal E}\times \MM$ valued Markov chains.
\end{theo}
\proof Since the process $(X(t), \pi(t))$ given by (\ref{9.10})-(\ref{9.10c}) is strongly Markov and the random times $S_k$ are $\mathcal{F}^X$-stopping times, the Markovity of $\xi_k$ follows from the definitions of $\XX_k$ and $\bar{\pi}_k$. Analogously, $\vt_k$ can also be shown Markov. \qed

\noi Let $M_k$ and $\bar{M}_k$ be the transition kernels of $\vt_k$ and $\xi_k$ respectively. By using the following measurable maps $\alpha(s,t, y^T_0):= t$ and $\varrho(s,t,y^T_0):= y^T_0(t)$ for $(s,t, y^T_0) \in {\cal E}$ and $s \le t$, we have $\alpha(\XX_k)= S_k$ and $\varrho(\XX_k)= X(S_k)$. Following the Remark \ref{9.rem2} we can express $\bar{\pi}_k$ in terms of $\XX_{k-1}, \XX_{k}$ and $\bar{\pi}_{k-1}$ in the following manner
\begin{eqnarray}\label{9.13}
\bar{\pi}_k= \tilde{\pi}^{\alpha(\XX_{k-1})}_{\alpha(\XX_k)}(\varrho(\XX_{k-1}), \bar{\pi}_{k-1}, \XX_k)
\end{eqnarray}
for $k=1, \ldots, n$
and
\begin{eqnarray}\label{9.13a}
\bar{\pi}^i_0 = \frac{d \eta(\cdot, i)}{d \sum_i \eta(\cdot, i)}(\XX_0)
\end{eqnarray}
for $i=1, \ldots, m$.

\begin{theo}\label{9.theo3}
(i) Let the conditional distribution $\bar{P}_k$ for $k= 1, \ldots, n$ be defined by
\begin{eqnarray}
\label{9.4a}\bar{P}_k \psi (e, \pi) := E [\psi(\XX_k) \mid \XX_{k-1}=e, \bar{\pi}_{k-1}= \pi]	
\end{eqnarray}
for all bounded continuous $\psi$, where $e \in {\cal E}$ and $\pi \in {\cal P}(\MM)$. Then for each $k=1,\ldots,n$
\begin{eqnarray*}\label{9.4}
\bar{M}_k \va (e, \pi) = \int \va(e', \tilde{\pi}^{\alpha(e)}_{\alpha(e')}(\varrho(e), \pi, e') ) \bar{P}_k(d e' \mid e, \pi)
\end{eqnarray*}
where $\tilde{\pi}^{\alpha(e)}_{\alpha(\cdot)}(\varrho(e), \pi, \cdot)$ is defined $\bar{P}_k(\cdot \mid e, \pi)$ almost surely.

(ii) $\bar{\pi}_k$ is measurable w.r.t $\{\XX_p: 0\le p\le k\}$, $k=0,\ldots, n$.
\end{theo}
\proof Since $\bar{M}_k$ is the transition kernel of $\xi_k$, using (\ref{9.13}) and (\ref{9.4a}), we have
\begin{eqnarray*}
% \nonumber to remove numbering (before each equation)
\lefteqn{\bar{M}_k \va (e, \pi)}\\
 &=& \int \va(e', \pi')\bar{M}_k(d e', d \pi' \mid e, \pi)\\
 &=& E [\va(\XX_k, \bar{\pi}_k) \mid \XX_{k-1}=e, \bar{\pi}_{k-1}= \pi]\\
 &=& E [E [\va(\XX_k, \bar{\pi}_k) \mid \XX_{k-1}, \XX_{k}, \bar{\pi}_{k-1}]\mid \XX_{k-1}=e, \bar{\pi}_{k-1}= \pi]\\
 &=& E [\va(\XX_k, \tilde{\pi}^{S_{k-1}}_{S_k}(\varrho(\XX_{k-1}), \pi_{k-1}, \XX_k)) \mid \XX_{k-1}=e, \bar{\pi}_{k-1}= \pi]\\
 &=& E [\va(\XX_k, \tilde{\pi}^{\alpha(e)}_{\alpha(\XX_k)}(\varrho(e), \pi, \XX_k)) \mid \XX_{k-1}=e, \bar{\pi}_{k-1}= \pi]\\
 &=& \int \va(e', \tilde{\pi}^{\alpha(e)}_{\alpha(e')}(\varrho(e), \pi, e') ) \bar{P}_k(d e' \mid e, \pi).
\end{eqnarray*}
And (ii) follows by (\ref{9.13a}) and using (\ref{9.13}) repeatedly to write $\bar{\pi}_k$ in terms of $\XX_k, \XX_{k-1}, \ldots, \XX_0$.\qed

\noi Let $ \mathcal{U}^j_k f (e, \theta) := E[f (\XX_j)\mid \XX_k =e, \bar{\theta}_k =\theta]$ and $ \bar{\mathcal{U}}^j_k f (e, \pi) := E[f (\XX_j)\mid \XX_k =e, \bar{\pi}_k= \pi]$ where $f$ be a bounded measurable function on ${\cal E}$, $j \geq k \geq 0$. We have the following lemma which plays the key role to prove the main result in Theorem \ref{9.16}.
\begin{lem}\label{9.lem1} Let $A$ be an event measurable w.r.t. $\XX_0, \ldots, \XX_{k-1}, \XX_k$, i.e., $A \in \sigma(\XX_0, \ldots, \XX_{k-1}, \XX_k)$. We have \\
(i)$E[\mathcal{U}^j_k f (\XX_k, \bar{\theta}_k)\mid \XX_k, \XX_{k-1}, \ldots, \XX_0]= E[\bar{\mathcal{U}}^j_k f (\XX_k, \bar{\pi}_k)\mid \XX_k, \XX_{k-1}, \ldots, \XX_0] $\\
(ii)$Var[\mathcal{U}^j_k f (\XX_k, \bar{\theta}_k)\mid I_A] \geq Var[\bar{\mathcal{U}}^j_k f (\XX_k, \bar{\pi}_k)\mid I_A] $.
\end{lem}
\proof (i) The result follows from the definition of $\mathcal{U}^j_k $ and $\bar{\mathcal{U}}^j_k $.

\noi (ii) By employing the rules of analysis of variance we have
\begin{eqnarray}
\no Var[\mathcal{U}^j_k f (\XX_k, \bar{\theta}_k)\mid I_A] &=& E[Var[\mathcal{U}^j_k f (\XX_k, \bar{\theta}_k) \mid \XX_k, \XX_{k-1}, \ldots, \XX_0, I_A]\mid I_A]\\
\label{9.11a} &&+ Var[E[\mathcal{U}^j_k f (\XX_k, \bar{\theta}_k) \mid \XX_k, \XX_{k-1}, \ldots, \XX_0, I_A]\mid I_A]\\
\no Var[\bar{\mathcal{U}}^j_k f (\XX_k, \bar{\pi}_k)\mid I_A] &=& E[Var[\bar{\mathcal{U}}^j_k f (\XX_k, \bar{\pi}_k)\mid \XX_k, \XX_{k-1}, \ldots, \XX_0, I_A]\mid I_A]\\
\label{9.11b} &&+ Var[E[\bar{\mathcal{U}}^j_k f (\XX_k, \bar{\pi}_k)\mid \XX_k, \XX_{k-1}, \ldots, \XX_0, I_A]\mid I_A].
\end{eqnarray}
Using (i) the second terms on the right in (\ref{9.11a}) and (\ref{9.11b}) are equal. The first term on the right in (\ref{9.11a}) is non-negative whereas that in (\ref{9.11b}) is zero (follows from Theorem \ref{9.theo3}). The result follows. ~~ \qed

\section{Feynman-Kac Distributions}
Define the potential functions
$$g_k(e) := I_{B_k}(\varrho(e))
$$
where $e \in {\cal E}$. Hence, the event $\{T_k \le T\}$ is equivalent to $\{ \prod_{p=0}^k g_p(\XX_p) =1\}$. Thus
\begin{eqnarray}\label{9.14}
I_{\{T_k \le T\}} = \prod_{p=0}^k g_p(\XX_p).	
\end{eqnarray}
Consider the Feynman-Kac distributions $\gamma_k$, the sequence of linear functionals on $L^\infty ({\cal E} \times \MM)$, given by
$$ \langle \gamma_k, \va \rangle: = E[\va(\XX_k, \bar{\theta}_k) \prod_{p=0}^k g_p(\XX_p)]$$
for $\va \in L^\infty ({\cal E} \times \MM)$. Since, $g_k$ are nonnegative, $\gamma_k$ are nonnegative bounded measures.
\begin{rem}
From Theorem \ref{9.theo2} we have
\begin{eqnarray*}
% \nonumber to remove numbering (before each equation)
\langle \gamma_k, \va \rangle &=& E\left[E[\va(\XX_k, \bar{\theta}_k) g_k(\XX_k)\mid \XX_{k-1}, \bar{\theta}_{k-1}]\prod_{p=0}^{k-1} g_p(\XX_p)\right]\\
  &=& E\left[M_k(\va g_k)(\XX_{k-1}, \bar{\theta}_{k-1})\prod_{p=0}^{k-1} g_p(\XX_p)\right]\\
  &=&\langle \gamma_{k-1}, M_k(\va g_k) \rangle\\
  &=&\langle \gamma_{k-1} M_k, \va g_k \rangle\\
 &=&\langle g_k \gamma_{k-1} M_k, \va \rangle.
 \end{eqnarray*}
Hence, $\gamma_k = g_k \gamma_{k-1} M_k$ or equivalently $\gamma_k= \gamma_{k-1}R_k$ where the nonnegative kernel $R_k$ is defined by $R_k(de',j\mid e, i):= M_k(de',j \mid e,i) g_k(e')$ for $k=1, \ldots, n$ and we also have $\gamma_0 = g_0\eta $.
\end{rem}

We are also interested in another set of Feynman-Kac distributions $\Gamma_k$, the sequence of linear functionals on $L^\infty ({\cal E} \times {\cal P}(\MM))$, given by
$$ \langle \Gamma_k, F \rangle: = E[F(\XX_k, \bar{\pi}_k) \prod_{p=0}^k g_p(\XX_p)]$$
for $F \in L^\infty ({\cal E} \times {\cal P}(\MM))$. Hence, analogously, $\Gamma_k$ are nonnegative bounded measures and $\Gamma_k = g_k \Gamma_{k-1} \bar{M}_k$ for $k=1, \ldots, n$ and $\Gamma_0 = g_0\eta_0 $ where $\eta_0$ is the distribution of $(X(0), \pi(0))$ given by \eqref{9.10b}-\eqref{9.10c}. Again for $\bar{R}_k(de',d \pi'\mid e, \pi):= \bar{M}_k(de',d\pi' \mid e,\pi) g_k(e')$, we have, $\Gamma_k= \Gamma_{k-1}\bar{R}_k$, $k=1, \ldots, n$. %The above derivations of recursion formulas for $\{\gamma_k: k=0, \ldots, n\}$ and $\{\Gamma_k: k=0, \ldots, n\}$ would be used to study the convergence of particle approximation of $P(\{T_n \le T\})= \langle \gamma_n, 1 \rangle = \langle \Gamma_n, 1 \rangle$ in the next section.

We note that the Feynman-Kac distributions $\gamma_k$ and $\Gamma_k$ are unnormalized distributions. The corresponding normalized distributions are $$\mu_k:= \frac{1}{\langle\gamma_k,1\rangle}\gamma_k
\txt{ ~and }
\bar{\mu}_k:= \frac{1}{\langle\Gamma_k,1\rangle}\Gamma_k
$$
respectively. Using (\ref{9.14}), we have $ \langle \mu_k, \va\rangle = E[\va(\XX_k, \bar{\theta}_k) \mid T_k \le T]$ and $ \langle \bar{\mu}_k, F\rangle = E[F(\XX_k, \bar{\pi}_k) \mid T_k \le T]$.

\begin{rem}\label{9.rem3}
We also define for any $k=0,1,\ldots, n$ kernels $R_{k+1:n}$ and $\bar{R}_{k+1:n}$ by
\begin{eqnarray*}
  R_{k+1:n} \va(e,i)&:=&R_{k+1}\cdots R_{n}\va(e,i) = E[\va(\XX_n, \bar{\theta}_n) \prod_{p=k+1}^n g_p(\XX_p)\mid \XX_k =e, \bar{\theta}_k=i]\\
  \bar{R}_{k+1:n} F(e,\pi)&:=&\bar{R}_{k+1}\cdots \bar{R}_{n}\va(e,\pi) = E[F(\XX_n, \bar{\pi}_n) \prod_{p=k+1}^n g_p(\XX_p)\mid \XX_k =e, \bar{\pi}_k=\pi]
\end{eqnarray*}
with the convention $R_{n+1:n} \va(e,i) =\va(e,i)$ and $\bar{R}_{n+1:n} F(e, \pi) =F(e, \pi) $.
\end{rem}

\section{Particle Approximations and Asymptotic Variance}
In this section we consider a numerical simulation technique by adopting multilevel splitting approach to estimate the probability $P(X(t) \in B \txt{ for some }t \in [0,T])$. In the earlier sections we have shown that by multilevel splitting the above probability can also be rewritten as $P(T_n\le T)$, $ \langle \gamma_n, 1 \rangle$ or $ \langle \Gamma_n, 1 \rangle$. We would study the particle approximation of $ \langle \Gamma_n, 1 \rangle$. The related random variables $\{\xi_k, k=0,1,\ldots, n\}$ are simulated by discretizing the stochastic differential equation (\ref{9.10})-(\ref{9.10c}).

We first generate $\xi^i_0 \sim \eta_0 $ $iid$ for $i =1 , \ldots, N$. We define the empirical distribution $\mu^N_0:= \frac{1}{N} \sum_{i=1}^N \delta_{\xi^i_0}$. Hence $ \langle \mu_0^N , \bar{M}_1(de, d\pi)\rangle = \sum_{i=1}^N \omega_0^i \bar{M}_1(de, d\pi \mid\xi^i_0)$ where $\omega_0^i =\frac{1}{N}$ for $i=1, \ldots, N$.
%Therefore $g_0 \cdot \eta_0^N = \sum_{i=1}^N \omega_0^i \delta_{\xi^i_0} =: \mu_0^N$ where $\omega_0^i = %\frac{g_0(\xi^i_0)}{\sum_{i=1}^N g_0(\xi^i_0)}$.

Next we generate $\xi^i_1 \sim \langle \mu_0^N , \bar{M}_1(de, d\pi)\rangle $ $iid$ for $i =1 , \ldots, N$. Among all the $N$ number of samples, according to splitting approach, we allot zero importance to those which fail to hit $B_1$. If none succeeds, the algorithm is stopped, otherwise we define the weighted empirical distribution $\mu^N_1 := \sum_{i=1}^N \omega_1^i\delta_{\xi^i_1}$ where $\omega_1^i = \frac{g_1(\xi^i_1)} {\sum_{i=1}^N g_1(\xi^i_1)}$.

In the above manner, for $k=2,\ldots, n$ we continue to generate $\xi^i_k \sim \langle \mu_{k-1}^N , \bar{M}_k(de, d\pi)\rangle $ $iid$ for $i =1 , \ldots, N$ and construct $\mu^N_k := \sum_{i=1}^N \omega_k^i\delta_{\xi^i_k}$ where $\omega_k^i = \frac{g_k(\xi^i_k)} {\sum_{i=1}^N g_k(\xi^i_k)}$, provided at least one of $\{\omega_k^i, i =1, \ldots, N\}$ is nonzero.

Let us define recursively,
\begin{eqnarray}\label{9.15}
\no \Gamma^N_k:= g_k \eta^N_k \langle \Gamma^N_{k-1}, 1\rangle \txt{ ~ with } \Gamma^N_0 := g_0 \eta^N_0
\end{eqnarray}
where $\eta^N_k := \frac{1}{N} \sum_{i=1}^N \delta_{\xi^i_k}$ for $k=0,1,\ldots, n$ and $\xi^i_k$ are generated as above.

\begin{theo}

The following holds
$$E|\frac{\langle \Gamma^N_{n}, 1\rangle}{\langle \Gamma_{n}, 1\rangle}-1| \le z^N_n \txt{ ~and } \sup_{F: \|F\|=1}E|\langle \mu^N_{n}-\bar{\mu}_n, F\rangle| \le 2 z^N_n
$$
where the sequence $\{ z^N_k\}$ satisfies the linear recursion
$$ z^N_k\le r_k (1+\frac{1}{\sqrt{N}}) z^N_{k-1}\frac{r_k}{\sqrt{N}} \txt{ ~ and } z^N_0 = \frac{r_0}{\sqrt{N}}
$$
for some constants $r_0, r_1, \ldots, r_n$.
\end{theo}
\proof The proof is analogous to that in \cite{CER} and \cite{CER1}. Therefore we omit the details. \qed

Hence, $\lim_{N \to \infty} E|\frac{\langle \Gamma^N_{n}, 1\rangle}{\langle \Gamma_{n}, 1\rangle}-1| = 0$. This shows that $\langle \Gamma^N_{n}, 1\rangle$ is an unbiased estimator of $\langle \Gamma_{n}, 1\rangle$. In particular we have the following result regarding the asymptotic variance. Since an analogous result is proved in \cite{CER} and \cite{CER1}, we state the result without a proof.
\begin{theo} Let ${\cal N} (a,b)$ denote the Normal distribution with mean $a$ and variance $b$. We have
$$ \lim_{N\to \infty}\sqrt{N}\left(\frac{\langle \Gamma^N_{n}, 1\rangle}{\langle \Gamma_{n}, 1\rangle}-1\right) \sim {\cal N} (0, \bar{V}_n) \txt{ ~and } \lim_{N\to \infty}\sqrt{N} \langle \mu^N_{n}-\bar{\mu}_n, F\rangle \sim {\cal N} (0, v_n(F))
$$
with
$$\bar{V}_n := \sum_{k=0}^n\left( \frac{\langle\eta_k , (g_k \bar{R}_{k+1:n}1)^2\rangle}{\langle \eta_k , g_k \bar{R}_{k+1:n}1\rangle^2}-1\right)
\txt{ ~and }
v_n(F) := \sum_{k=0}^n \frac{\langle\eta_k , |g_k \bar{R}_{k+1:n}(F - \langle \bar{\mu}_n, F\rangle)|^2 \rangle}{\langle\eta_k , g_k \bar{R}_{k+1:n}1\rangle^2}
$$
where $\eta_k := \bar{\mu}_{k-1}\bar{M}_k$ for $k=1, \ldots, n$.
\end{theo}

Since, $g_k$ is an indicator function, ${g_k}^2= g_k$ for each $k$. Using this, we have the following simplification of $\bar{V}_n$
\begin{eqnarray}\label{9.17}
\bar{V}_n = \sum_{k=0}^n \left(\frac{1}{p_k}-1\right) + \sum_{k=0}^n\frac{1}{p_k} \frac{Var (\bar{R}_{k+1:n}1, \bar{\mu}_k)}{\langle \bar{\mu}_k ,\bar{R}_{k+1:n}1\rangle^2}
\end{eqnarray}
with $p_k:= \langle \eta_k, g_k\rangle = P(T_k \le T \mid T_{k-1} \le T ), k=1, \ldots, n$. We also define
\begin{eqnarray}\label{9.17a}
V_n := \sum_{k=0}^n \left(\frac{1}{p_k}-1\right) + \sum_{k=0}^n\frac{1}{p_k} \frac{Var (R_{k+1:n}1, \mu_k)}{\langle \mu_k ,R_{k+1:n}1\rangle^2}
\end{eqnarray}
analogous to (\ref{9.17}). Therefore, $V_n$ is the asymptotic variance of a similar particle approximation of $\langle \gamma_{n}, 1\rangle$ where the random variables $\{\vartheta_k, k=0,1,\ldots, n\}$ are simulated by discretizing the stochastic differential equation (\ref{9.1})-(\ref{9.1b}). The details of this alternative particle approximation procedure is analogous to that, stated in this section and therefore, we skip the details. We wish to compare both of the particle approximations by comparing the corresponding asymptotic variances.

\begin{theo} \label{9.16} $\bar{V}_n \leq 	 V_n.$
\end{theo}
\proof The expressions \eqref{9.17} and \eqref{9.17a} are being used for the comparison. From Remark \ref{9.rem3}, we have
$$Var (\bar{R}_{k+1:n}1, \bar{\mu}_k) = Var[P[T_n \le T \mid \XX_k, \bar{\pi}_k]\mid T_k \leq T]
$$
and
$$Var (R_{k+1:n}1, \mu_k) = Var[P[T_n \le T \mid \XX_k, \bar{\theta}_k]\mid T_k \leq T].
$$
We rewrite $P[T_n \le T \mid \XX_k, \bar{\pi}_k] = E[I_{[0,T]} \circ \alpha (\XX_n)\mid \XX_k, \bar{\pi}_k] = \bar{\mathcal{U}}^n_k f (\XX_k, \bar{\pi}_k)$, where $f =I_{[0,T]} \circ \alpha $, a bounded measurable map. Similarly, we rewrite $P[T_n \le T \mid \XX_k, \bar{\theta}_k] = \mathcal{U}^n_k f (\XX_k, \bar{\theta}_k)$
Thus by using Lemma \ref{9.lem1} we get $Var (\bar{R}_{k+1:n}1, \bar{\mu}_k) \leq Var (R_{k+1:n}1, \mu_k)$ and $ \langle \bar{\mu}_k ,\bar{R}_{k+1:n}1\rangle =\langle \mu_k ,R_{k+1:n}1\rangle$. Hence the result. \qed

\begin{rem}\label{9.rem4}
The inequality in Theorem \ref{9.16} establishes the fact that the use of Wonham filter reduces the asymptotic variance of the particle approximation. In this section we have selected a standard particle approximation scheme among many others. There are, indeed, some other schemes leading a lesser asymptotic variance \cite{DEL}. But we wish to emphasize that, the use of Wonham filter does not restrict one from using those schemes. In fact the use of Wonham filter enables further reduction in asymptotic variance for any other choice of particle approximation scheme. Here we illustrate this fact by considering an alternative scheme.

We first generate initial particles $\xi^i_0 \sim \eta_0$ $iid$ for $i =1 , \ldots, N$. We define the empirical distributions $\hat{\mu}^N_0 = \hat{\eta}^N_0 := \frac{1}{N} \sum_{i=1}^N \delta_{\xi^i_0}$. Hence $ \langle \hat{\mu}_0^N , \bar{M}_1(de, d\pi)\rangle = \frac{1}{N} \sum_{i=1}^N \bar{M}_1(de, d\pi \mid \xi^i_0)$ for $i=1, \ldots, N$.
%Therefore $g_0 \cdot \eta_0^N = \sum_{i=1}^N \omega_0^i \delta_{\xi^i_0} =: \mu_0^N$ where $\omega_0^i = %\frac{g_0(\xi^i_0)}{\sum_{i=1}^N g_0(\xi^i_0)}$.

Next we generate $\xi^i_1 \sim \langle \hat{\mu}_0^N , \bar{M}_1(de, d\pi)\rangle $ $iid$ for $i =1 , \ldots, N$ and let $ \hat{\eta}^N_1:= \frac{1}{N} \sum_{i=1}^N \delta_{\xi^i_1}$. If none of the $N$ number of particles succeeds to hit $B_1$, the algorithm is stopped. Otherwise we resample the particle system $(\xi^1_1, \ldots, \xi^N_1)$ to obtain the \emph{selected} particle system $(\hat{\xi}^1_1, \ldots, \hat{\xi}^N_1)$, with $\hat{\xi}^i_1 \sim \bar{M}_1(\hat{\eta}^N_1, \xi^i_1, d\xi')$ where, $\bar{M}_1(\hat{\eta}, (e,\pi), d(e,\pi)') := g_1(e)\delta_{(e,\pi)}(d(e,\pi)') + (1-g_1(e))g_1\cdot\hat{\eta}(d(e,\pi)')$. The corresponding empirical distribution $\hat{\mu}^N_1 := \frac{1}{N}\sum_{i=1}^N \delta_{\hat{\xi}^i_1}$.

In the above manner, for $k=2,\ldots, n$ we continue to generate $\xi^i_k \sim \langle \hat{\mu}_{k-1}^N , \bar{M}_k(de, d\pi)\rangle $ $iid$ for $i =1 , \ldots, N$ and define $\hat{\eta}^N_k := \frac{1}{N}\sum_{i=1}^N \delta_{\xi^i_k}$, $\bar{M}_k(\hat{\eta}, (e,\pi), d(e,\pi)') := g_k(e)\delta_{(e,\pi)}(d(e,\pi)') + (1-g_k(e))g_k\cdot\hat{\eta}(d(e,\pi)')$. Next the particles $\hat{\xi}^i_k \sim \bar{M}_k(\hat{\eta}^N_k, \xi^i_k, d\xi')$ are produced for $i=1, \ldots, N$; provided at least one of $\{g_k(\xi_k^i), i =1, \ldots, N\}$ is nonzero. The $k$th empirical distribution $\hat{\mu}^N_k := \frac{1}{N}\sum_{i=1}^N \delta_{\hat{\xi}^i_k}$ is also obtained.

Let us define recursively for $k=1,\ldots, n$
\begin{eqnarray}\label{9.19}
\no \hat{\Gamma}^N_k:= g_k \hat{\eta}^N_k \langle \hat{\Gamma}^N_{k-1}, 1\rangle \txt{ ~ with } \hat{\Gamma}^N_0 := g_0 \hat{\eta}^N_0.
\end{eqnarray}
Then $\langle \hat{\Gamma}^N_{n}, 1\rangle$ is also an unbiased estimator of $\langle \Gamma_{n}, 1\rangle$. %(see \cite{?}).
The corresponding asymptotic variance $\bar{\hat{V}}_n$ can be expressed as
\begin{eqnarray}\label{9.20}
\bar{\hat{V}}_n = \sum_{k=0}^n \left(\frac{1}{p_k}-1\right) + \sum_{k=0}^n\frac{1}{p_k} \frac{Var (\bar{R}_{k+1:n}1, \bar{\mu}_k)}{\langle \bar{\mu}_k ,\bar{R}_{k+1:n}1\rangle^2}(1-{p_k}^2)
\end{eqnarray}
where $p_k$ are as before. We also define analogously
\begin{eqnarray}\label{9.20a}
\hat{V}_n := \sum_{k=0}^n \left(\frac{1}{p_k}-1\right) + \sum_{k=0}^n\frac{1}{p_k} \frac{Var (R_{k+1:n}1, \mu_k)}{\langle \mu_k ,R_{k+1:n}1\rangle^2}(1-{p_k}^2)
\end{eqnarray}
Therefore, $\hat{V}_n$ is the asymptotic variance of a similar particle approximation of $\langle \gamma_{n}, 1\rangle$ where the random variables $\{\vartheta_k, k=0,1,\ldots, n\}$ are simulated by discretizing the stochastic differential equation (\ref{9.1})-(\ref{9.1b}). As before, using Lemma \ref{9.lem1} we compare between (\ref{9.20}) and (\ref{9.20a}) to get $\bar{\hat{V}}_n \le \hat{V}_n$.
\end{rem}

\section{Conclusion}
Though the use of Wonham filter causes variance reduction but it requires more computational steps for numerical simulation. Therefore to be sure of an overall advantage of this method, it is important to know the difference between variances. Since the difference or the amount of reduction depends on the parameters of the underlying hybrid process, a general study is quite involved. Nevertheless we can check for some extreme cases. An example is given below.

\begin{ex}
(i) Let $\Lambda(x)$ be a constant null matrix. We also assume that $\eta$ be a product measure of the form $\eta(dx,i)= \nu(dx) \delta_{j}(i)$ then $Var[ f (\XX_k, \bar{\theta}_k) \mid \XX_k, \XX_{k-1}, \ldots, \XX_0] =0$ for all admissible $k$ and $f$. Hence, $\bar{V}_n =V_n$. On the other hand, both the drift and the noise coefficients in (\ref{9.10a}) vanish for $t \ge 0$ and results in $\pi(t)$ a constant vector equal to $1_j$. Therefore, for such trivial case, the use of Wonham filter neither reduces variance nor enhances the computational task.
\end{ex}

%\section{Numerical Examples}

\end{document}